\documentclass[12pt]{article}
\usepackage{amsmath,amsthm,amssymb,amsfonts,fullpage}

\usepackage{amssymb,amsmath,latexsym, bbm}
\usepackage{epsfig}
\newtheorem{thm}{Theorem}[section]

\newtheorem{lem}[thm]{Lemma}
\newtheorem{propn}[thm]{Proposition}
\newtheorem{rem}[thm]{Remark}\newtheorem{assump}[thm]{Assumption}

\def\qed{{\hfill $\Box$ \bigskip}}

\def\bR{\mathbb{R}}

\def\G{\mathcal{G}}

\def\lam{\lambda}

\def\1{\mathbbm{1}}
\def\bar{\overline}
\def\eps{\varepsilon}

\def\<{\langle}
\def\>{\rangle}
\def\bee{\begin{equation}}
  \def\eee{\end{equation}}
  \makeatletter\@addtoreset{equation}{section}
  
  \makeatother
\begin{document}

  \bibliographystyle{plain}
\title{\Large \bf 
Fluctuations of recentered maxima of discrete Gaussian Free Fields on 
a class of recurrent graphs}
\author{{\bf Takashi~Kumagai}\thanks{Research partially supported by the
Grant-in-Aid for Scientific Research (B) 22340017.} 
\and   {\bf Ofer Zeitouni}
\thanks{Research partially supported by NSF grant \#DMS-1106627 and
a grant from the Israel Science Foundation.}}
\date{February 8, 2013}
\maketitle
\begin{abstract} We provide conditions that ensure that the recentered
  maximum of
the Gaussian free field on a sequence of graphs fluctuates at the same order
as the field at the point of maximal variance. In particular, on a sequence of
such graphs the recentered
maximum is not tight, similarly to the situation in
$\mathbb{Z}$ but in contrast with the situation in $\mathbb{Z}^2$.
We show that
our conditions cover 
a large class of ``fractal'' graphs.
\end{abstract}
%

  \section{Introduction}  
  The study of the maxima of Gaussian fields has a rich
  history, which we will not attempt to survey here. The general
  theory was developed in the 70s and 80s, and an excellent
  account can be found in \cite{LT}. However, general results concerning the 
  order of fluctuations of the maximum are lacking.

  In recent years, a special effort has been directed toward the study of the
  so called Gaussian free field (GFF) on various graphs. While we postpone
  the general definition to the next section, we discuss in this introduction the
  special case of the GFF on subsets $V_N=([-N,N]\cap \mathbb{Z})^d$, with
  Dirichlet boundary conditions. These are random fields $\{X_x\}_{x\in V_N}$
  indexed by points
  in $V_N$, with joint density (with respect to Lebesgue measure)
  proportional to
  $$ \exp\left(- c\sum_{x\sim y} (X_x-X_y)^2\right)\,,$$
  with the sum  over neighbors in $V_N$, and $X_x=0$ for $x\in \partial V_N$.
  (An alternative description involving the Green function of random walk
  on $V_N$ is given below in Section \ref{sect2}; see also \cite{She}
  for a very readable introduction to GFFs in a continuous setting.)
  With $X_{N,d}^*$ denoting the maximum of the GFF on $V_N$ in dimension $d$, 
  it is not hard to see that
  $X_{N,d}^*$ is of order $\sqrt{N}$ for $d=1$, order $\log N$
  for $d=2$, and order $(\log N)^{1/2}$ for $d\geq 3$. Moreover, a consequence
  of the Borell-Tsirelson inequality (see \cite{LT}) is that for $d\geq 3$,
  since simple random walk is transient on $\mathbb{Z}^d$,
  the fluctuations of $X_{N,d}^*$ are at most of order $1$, while
  for $d=1$ the fluctuations of $X_{N,1}^*$ are of the same order as
  $X_{N,1}^*$, i.e. of order $\sqrt{N}$. The critical case $d=2$ was
  settled only recently \cite{BZ10}, where it was shown that the fluctuations
  of $X_{N,2}^*$ are also of order $1$. This raises naturally the question
  of determining for which sequences of graphs is the sequence of
  recentered maxima of the GFF tight.

Our goal in this paper is to exhibit a class of sequences of
graphs, which are fractal-like and  for which the maximum
of the GFF fluctuates
at the same order as the maximum itself, and both are
of the order of the 
maximal standard deviation of the GFF
in the graph. In that respect, the behavior
of the maximum 
is similar to that of $X_{N,1}^*$.
For this class of graphs, we also show that the \textit{cover time} of 
the graph, measured in terms of the (square root of the)
local time at a fixed vertex, also does not concentrate.
(We note in passing that for the cover time of $V_N$ in two dimensions, it 
is, to the best of our knowledge, an open problem to decide whether this
quantity concentrates or not.)

The structure of the paper is as follows. In the next section, we introduce the GFF on general graphs and state Assumption \ref{thm:assump1}
that characterizes the graphs 
which we investigate; the main feature is a relation between the graph distance
and the resistance, and  control of the covering number 
of the graph in terms of resistance distance. 
We then state our main result,
Theorem \ref{thm:main thm}, concerning fluctuations of the maximum of the GFF.
We also state Proposition
    \ref{thm:main prop} concerning the cover time of the graphs. Proofs of
    the theorem and proposition are given in Section
    \ref{sect:proofs}. The heart of the paper is then Section 
    \ref{sect:examples}, where we show that certain naturally constructed 
    fractal-like graphs satisfy our assumptions. In particular, this is the
    case for the standard Sierpinski carpets in two dimensions
    and gaskets in all dimensions.

    \textbf{Notation}
    Throughout the paper, we use $c_1, c_2, \cdots$ to denote generic
  constants, independent of $N$,
  whose exact values are not important and may  change from
  line to line. We write $a_n\asymp b_n$ if there exist constants $c_1, c_2>0$ 
  such that
  $c_1b_n\le a_n\le c_2b_n$ for all $n\in \mathbb N$. 
  \section{Framework}\label{sect2}
  We first introduce general notation for finite graphs with a \lq wired' boundary and their 
  associated resistance. Let $G=(V(G),E(G))$  be a connected 
  (undirected) finite graph with at least two vertices, 
  where $V(G)$ denotes the vertex set and $E(G)$ the edge set of $G$. 
  Let $d_G$ be the graph distance, that is, 
  $d_G(x,y)$ is the number of edges in the shortest path from $x$ to 
  $y$ in $G$. Define a symmetric weight function 
  $\mu^G:V(G)\times V(G)\rightarrow \mathbb{R}_+$ that satisfies 
  $\mu^G_{xy}>0$ if and only if $\{x,y\}\in E(G)$. For 
  $B\subset G$ with $B\ne G$ and for distinct $x,y\in V(G)$ not both in $B$, 
  we define the 
  \textit{resistance} between $x$ and $y$ 
  by
  \[R_B(x,y)^{-1}:=\inf\{\,\frac 12\sum_{w,z\in V(G)} (f(w)-f(z))^2\mu_{wz}^G: 
  f(x)=1, f(y)=0, f|_B= \mbox{\rm constant}\}.\]
  We set $R_B(x,x)=0$, $R_B(x,y)=0$ if $x,y\in B$
  and, for $x\in V(G)\setminus B$,
  we define $R_B(x,B)=R_B(x,y)$ for any $y\in B$. 
  We write $R(x,y):=R_{\emptyset}(x,y)$. 
  
  The resistance 
  $R_B(\cdot,\cdot)$
  is the  resistance of the following electrical 
  network with a \lq wired' boundary:  
    Consider the graph $\bar G$ obtained by combining all vertices in $B$ to a
single vertex $b$, that is $V(\bar G)=(V(G)\setminus B)\cup \{b\}$ and
\begin{eqnarray*}
  E(\bar G)&=&\{\{x,y\}: \{x,y\}\in E(G), x,y\in G\setminus B\}
\\&&
\bigcup \{\{x,b\}: x\in G\setminus B,\, \exists y\in B \,\mbox{\rm 
with} \, \{x,y\}\in E(G)\}\,.
\end{eqnarray*}
Define the modified symmetric weight function
$$ \mu^{\bar G}_{xy}=\left\{
\begin{array}{ll}
  \mu^{G}_{xy},& 
  x\in V(\bar G)\setminus \{b\}, \,y\in V(\bar G)\setminus \{b\},\\
  \sum_{z\in B} \mu^{G}_{xz},&
   x\in V(\bar G)\setminus \{b\}, \,y=b\,,
\end{array}
\right.$$
and set as before $\mu^{\bar G}_x=\sum_{y\in V(\bar G)} 
\mu^{\bar G}_{xy}$. 
Let $\{\bar w_t\}_{t\geq 0}$ be the continuous time random walk on $\bar G$ 
such that the holding time at a vertex is $\exp(1)$, and the jump probability is given by 
$\mu_{x,y}/\mu_x$. Let $$L_t^{x,N}=\frac{1}{\mu^{\bar G}_x}
\int_0^t {\bf 1}_{\{\bar w_s=x\}}ds$$
denote the (weight normalized) local time at $x$. 
  
  Now, let $\{G^N\}_{N\ge 1}$ be a sequence of finite connected graphs 
  such that $|G^N|\ge 2$ for all $N\ge 1$ and 
  $\lim_{N\to\infty}|G^N|=\infty$. For each $G^N=(V(G^N), E(G^N))$, 
 we take a symmetric weight function $\mu^{G^N}$, a boundary 
 $B^N\subset G^N$ with $B^N\ne G^N$, 
 and the corresponding continuous time Markov chain $\{\bar w_t^N\}_{t\ge 0}$ 
with the wired boundary condition on $B^N$  as above. 
We assume that $G^N\setminus B^N$ is connected. 
Let $T^N:=\min\{t\ge 0: \bar w^{N}_t= b\}$, 
and define, for each 
$x,y\in V(G^N)\setminus B^N$,
$\G_N(x,y)=E_{G^N}^x[L_{T^N}^{y,N}]$ 
where
$E_{G^N}^x$ denotes the expectation with respect to $\bar w_t^{N}$ started at $x$.
For $z\in B^N$, we set $X_z^N\equiv 0$. 
The \textit{Gaussian free field} (GFF for short) on $G^N$ (with boundary $B^N$)
is the zero-mean Gaussian field $\{X_z^N\}_{z\in V(G^N)}$ 
with covariance $\G_N(\cdot,\cdot)$. 
It can be easily checked (using for instance
\cite[Lemma 2.1]{DLP}, \cite[Proposition 3.6]{TKSF}) that 
\[E[(X_x^N-X_y^N)^2]=R_{B^N}(x,y).\]
Let $h :\mathbb N\to \mathbb N$ be a strictly
increasing function with $h(0)=0$, 
that satisfies the following doubling property: there exist 
$0<\beta _{1}\leq \beta _{2}<\infty$ and $C>0$ such that, for 
all $0<r\leq R<\infty $, 
\begin{equation}{C^{-1}}\left( \frac{R}{r}\right) ^{\beta _{1}}\leq 
  \frac{h(R)}{h (r)}\leq C\left( \frac{R}{r}\right) ^{\beta _{2}}.  
  \label{eq:upvol}\end{equation}
 
  We assume the following. 
 \begin{assump}\label{thm:assump1} 
   There exist $\alpha>0$ and $c_1,c_2,c_3>0$ such that the following 
   hold for all large $N$.\\
   {\rm (i) } $R_{B^N}(x,y)\le c_1h(d_{G^N}(x,y))$ for all $x,y\in G^N$. \\
   {\rm (ii) } $\max_{x\in G^N} R_{B^N}(x,B^N)\ge 
   c_2\max_{x\in  G^N}h(d_{G^N}(x,B^N))$  
   for all $x\in  G^N$.\\
   {\rm (iii) } ${\mathcal N}_{G^N}(\delta d_{max}^N)\le 
   c_3\delta^{-\alpha}$ 
   for all $\delta \in (0,1]$ where 
   $d_{max}^N:=\max_{x\in  G^N}d_{G^N}(x,B^N)$ and 
   ${\mathcal N}_{G^N}(\varepsilon)$ is the minimal number of 
   $d_{G^N}$-balls of radius $\varepsilon$ needed to cover $G^N$. 
   Furthermore, $d_{max}^N\to\infty$ as $N\to\infty$.    \end{assump}
   Let $X_N^*=\max_{z\in V(G^N)}X_z^{N}$ and define 
   $\tilde X_N=X_N^*/\bar \sigma_N$, where 
   $\bar \sigma_N=(\max_{z\in  G^N}E[(X_z^N)^2])^{1/2}$. 
   Note that $\bar \sigma_N^2=\max_{x\in G^N} R_{B^N}(x,B^N)$, and $\lim_{N\to\infty}\bar \sigma_N=\infty$
  under Assumption \ref{thm:assump1}(iii).

   \begin{thm}\label{thm:main thm}
     Under Assumption \ref{thm:assump1}, there exist constants $A,B,A'>0$ and a 
     function $g: (0,\infty)\to (0,1)$ such that the following 
     holds for all $N$ large. 
     \bee\label{eq:nolln}
     P(\tilde X_N<A)>B,~~~ P(\tilde X_N>c)\ge g(c)~~\forall c>0,~~~
     E(\tilde X_N)\leq A'.\eee
   \end{thm}
   In particular, under Assumption \ref{thm:assump1}, 
      $\{X_N^*-EX_N^*\}_N$ fluctuates with order $\bar \sigma_N$ and  
     therefore  it is not tight. 

     \begin{rem} We
       stated Assumption \ref{thm:assump1} with respect to the graph
  distance in $G^N$, because this will be easiest to check in the applications.
  However, one should note that the proof of Theorem \ref{thm:assump1} does not
  depend on the particular metric chosen, as long as the metric satisfies
  the assumption.
  In particular, if we choose $R_{B^N}(\cdot,\cdot)$ as the metric, Assumption \ref{thm:assump1}
  (i), (ii) turns out to be trivial with $h(s)=s$, and the assumption boils down to 
  ${\mathcal N}_{R_{B^N}}(\delta \bar \sigma_N^2)\le c_3\delta^{-\alpha}$ 
   for all $\delta \in (0,1]$ and $\lim_{N\to\infty}\bar 
   \sigma_N=\infty$, where 
    ${\mathcal N}_{R_{B^N}}(\varepsilon)$ is the minimal number of 
   $R_{B^N}$-balls of radius $\varepsilon$ needed to cover $G^N$. 
\end{rem}
   
     In a recent seminal work,
     \cite{DLP} have established 
     a close relation between the expectation of
     the maximum of the GFF on general graphs
and the expected cover time of these graphs by random walk.
Under the assumptions of Theorem \ref{thm:main thm}, one can also
derive information on the fluctuations of the cover time, as follows.
Define the 
\textit{cover time} of $\bar G^N$ as
$$\tau_{\mbox{\rm cov}}^N=\inf\{t>0: L_t^{x,N}>0,\, \forall x\in \bar G^N\}\,.$$
It is easy to see
$$\tau_{\mbox{\rm cov}}^N=\inf\{t>0: \forall x\in \bar G^N, \exists s\le t \,\mbox{ such that }\, 
\bar w_s^N=x\}. $$
We will consider the square-root of the
normalized local time at $B^N$ at cover time, i.e.
the random variable $L^N:=\sqrt{L_{\tau_{\mbox{\rm cov}}^N}^{b, N}}$. One expects
(see \cite{DLP}) that $L^N$ should behave similarly to $|X_N^*|$.
In the special case of $G^N$ being the rooted at $b$ binary tree of depth $N$,
this was confirmed in \cite{BZrec}.
In our setup here, this is confirmed in the following proposition. 
\begin{propn}\label{thm:main prop}
  With notation as above and under Assumption
  \ref{thm:assump1}, the conclusion of Theorem
  \ref{thm:main thm} hold with $L^N/\bar \sigma_N$ replacing $\tilde X_N$.
  \end{propn}
    \section{Proofs of Theorem \ref{thm:main thm} and Proposition
    \ref{thm:main prop}}
    \label{sect:proofs}
    We begin with the proof of Theorem \ref{thm:main thm}.\\
    {\it Proof of Theorem \ref{thm:main thm}}: Let 
    $\tilde d(x,y)=(E[(X_x^N-X_y^N)^2])^{1/2}/\bar \sigma_N
=
    R_{B^N}(x,y)^{1/2}/\bar \sigma_N$. 
    Then, using Assumption \ref{thm:assump1} (i),(ii), 
    there exists
$c>0$ such that for all $x,y\in  G^N$ with 
    $d_{G^N}(x,y)\le d_{max}^N$ and all $N\in \mathbb N$, 
    \[\tilde d(x,y)=\Big(\frac{R_{B^N}(x,y)}{\bar \sigma_N^2}\Big)^{1/2}\le 
    c
\Big(\frac{h(d_{G^N}(x,y))}{h(d_{max}^N)}\Big)^{1/2}
\le 
    cC\Big(\frac{d_{G^N}(x,y)}{d_{max}^N}\Big)^{\beta_1/2}.
\]
    Thus, denoting ${\mathcal N}_{\tilde d}(\varepsilon)$ the minimal number 
    of $\tilde d$-balls of radius $\varepsilon$ needed to cover $G^N$, 
    we have
    \[{\mathcal N}_{\tilde d}(cC\delta^{\beta_1/2})\le 
    {\mathcal N}_{ G^N}(\delta d_{max}^N)\le c_3\delta^{-\alpha},\]
    where we used  Assumption \ref{thm:assump1} (iii) in the second inequality. Rewriting this, we have
    ${\mathcal N}_{\tilde d}(\varepsilon)\le c'\varepsilon^{-2\alpha/\beta_1}$, where $c'>0$ is independent of $N$. 
    Set $\gamma=2\alpha/\beta_1$. We can apply \cite[Theorem 5.2]{Adl} to 
    deduce that there exist $\lambda_0>0$ and $N_0$ such that for all $\lambda>\lambda_0$,
    $\eps>0$ and $N>N_0$,
    \[
    P(\tilde X_N>\lambda)\le C_\gamma
    \lambda^{\gamma+1+\eps}\Psi(\lambda),    
     \]
   where $C_\gamma\ge 1$ does not depend on $N$ and 
   $\Psi(\lambda)=(2\pi)^{-1/2}\int_\lambda^\infty e^{-x^2/2}dx$. On the 
   other hand, let $x^*_N$ be such that $E(X_{x^*_N}^2)=\bar \sigma_N^2$.
   Then, for any $\lambda>0$,
   $$P(\tilde X_N>\lambda)\ge P(X_{x_N^*}^N>\lambda \bar \sigma_N)=
    \Psi(\lambda)\,.$$
   The estimates in 
   (\ref{eq:nolln}) are easy consequences of the last two displayed
   inequalities. \qed 
   
   \noindent We turn to the analysis of cover times.\\
   {\it Proof of Proposition \ref{thm:main prop}}: 
   The upper bound in the proposition is a consequence of the 
   Eisenbaum-Kaspi-Marcus-Rosen-Shi isomorphism theorem \cite{EKMRS},
   as was observed in \cite{DLP}: indeed, by 
   \cite[Eq. (20),(21)]{DLP} and using the last
   estimate in \eqref{eq:nolln}, 
   there exist constants $c_1,c_2>0$ so that
   with $t=\theta \bar \sigma_N^2$, and all $\theta$ large enough,
     \begin{equation}
       \label{eq-oof2}
       P(\min_x L_{\tau^N(t)}^{x,N}\leq t/2)\leq c_1e^{-c_2 \theta}
     \end{equation}
    while
     \begin{equation}
      \label{eq-oof3}
      P(\max_x L_{\tau^N(t)}^{x,N}\geq 2t)\leq c_1e^{-c_2 \theta}\,,
    \end{equation}
    where $\tau^N(t):=\inf\{s> 0: L_{s}^{b,N}>t \}$. 

     On the event 
    $\{\min_x L_{\tau^N(t)}^{x,N}\geq t/2\}$ we have that 
    $\tau^N(t)\geq 
     \tau_{\mbox{\rm cov}}^N$.
    Thus, on the event 
    $$\{\min_x L_{\tau^N(t)}^{x,N}\geq t/2\}\cap 
    \{\max_x L_{\tau^N(t)}^{x,N}\leq 2t\}\,,$$
     one has that 
     \begin{equation}
       \label{eq-oof1}
       L_{\tau_{\mbox{\rm cov}}^N}^{b,N}\leq
     L_{\tau^N(t)}^{b,N}\leq\max_x
     L_{\tau^N(t)}^{x,N}\leq 2t\,.\end{equation}
     In particular, \eqref{eq-oof2}, \eqref{eq-oof3} and \eqref{eq-oof1}
     imply  that
     $EL^N/\bar \sigma_N$ is bounded uniformly.

     To estimate $L^N$ from below, we use the Markov property.
     Let $x^*\in V(\bar G^N)$ be such that $R_{B^N}(x^*,B^N)=\bar \sigma_N^2$
     and let $T_{x^*}=\inf\{t: \bar w_t^N=x^*\}$.
     Since $\tau_{\mbox{\rm cov}}^N\geq T_{x^*}$, we have that
     $L^N\geq \sqrt{L_{T_{x^*}}^{b,N}}$. We decompose the walk
     $\bar w_t^N$ according to excursions from $b$: the probability to hit
     $x^*$ during one excursion (see e.g. \cite[Ch. 2]{LP})
     is 
     $$p_N=\frac{1}{\bar \sigma_N^2 \mu_N}\,,$$
     where $\mu_N= \mu^{\bar G^N}_b\,.$
      Therefore,
     $${L_{T_{x^*}}^{b,N}}
     \stackrel{d}{=} \frac{1}{\mu_N}\sum_{i=1}^{Z_N} \mathcal{E}_i\,,$$
     where $Z_N$ is geometric of parameter $p_N$ and
     $\mathcal{E}_i$ are standard
     independent exponential random variables.
     Note that $E L_{T_{x^*}}^{b,N}=\bar \sigma_N^2$. 

     Consider now a parameter $\xi>0$. 
     We have that
     \begin{eqnarray*}
       P(
     L_{T_{x^*}}^{b,N}\geq \xi \bar\sigma_N^2)
     &\geq & P(Z_N\geq  \xi/p_N)
     P\left(\frac1{\mu_N}\sum_{i=1}^{\xi/p_N}\mathcal{E}_i>
     \xi \bar \sigma_N^2\right)\\
     &\geq &
     P(Z_N\geq \xi/p_N)P\left(\frac{p_N}{\xi}
     \sum_{i=1}^{\xi/p_N}\mathcal{E}_i\geq 1\right)=:P_1P_2\,.
     \end{eqnarray*}
Note that from the properties of the geometric distribution,
regardless of $p_N$ we have that $P_1\geq c_1(\xi)>0$.
 On the other hand, if $p_N\to 0$ then
     $p_N \sum_{i=1}^{1/p_N}\mathcal{E}_i\to 1$ a.s., and in any case
     we also have that $P_2\geq c_2(\xi)>0$. We conclude that
     $$P(L^N\geq \sqrt{\xi}\bar \sigma_N)\geq 
     c_1(\xi) c_2(\xi)\,.$$
     \qed

   \section{Examples}
      \subsection{Nested fractal graphs and strongly recurrent Sierpinski carpet graphs}\label{sssec}
Let $\{\psi_i\}_{i=1}^K$ be
a family of $L$-similitudes on $\mathbb{R}^d$ for some $L>1$, that is, for each $i$, $\psi_i$ is a map from $\mathbb{R}^d$ to $\mathbb{R}^d$ such that 
$\psi_i(x)=L^{-1}U_i x +\gamma_i,~ x\in {\mathbb R}^d$, where $U_i$ is a unitary map and $\gamma_i\in {\mathbb R}^d$. 
We assume that $\{\psi_i\}_{i=1}^K$ satisfies the open set condition,
namely there exists a non-empty bounded set $O\subset\mathbb{R}^d$ such that $\{\psi_i(O)\}_{i=1}^K$ are disjoint and $\cup_{i=1}^K\psi_i(O)\subset O$. Since $\{\psi_i\}_{i=1}^K$ is a family of contraction maps, there exists a unique non-empty compact set $F$ such that $F=\cup_{i=1}^K \psi_i(F)$. 
We assume that $F$ is connected.

   \label{sect:examples}
   \centerline{\epsfig{file=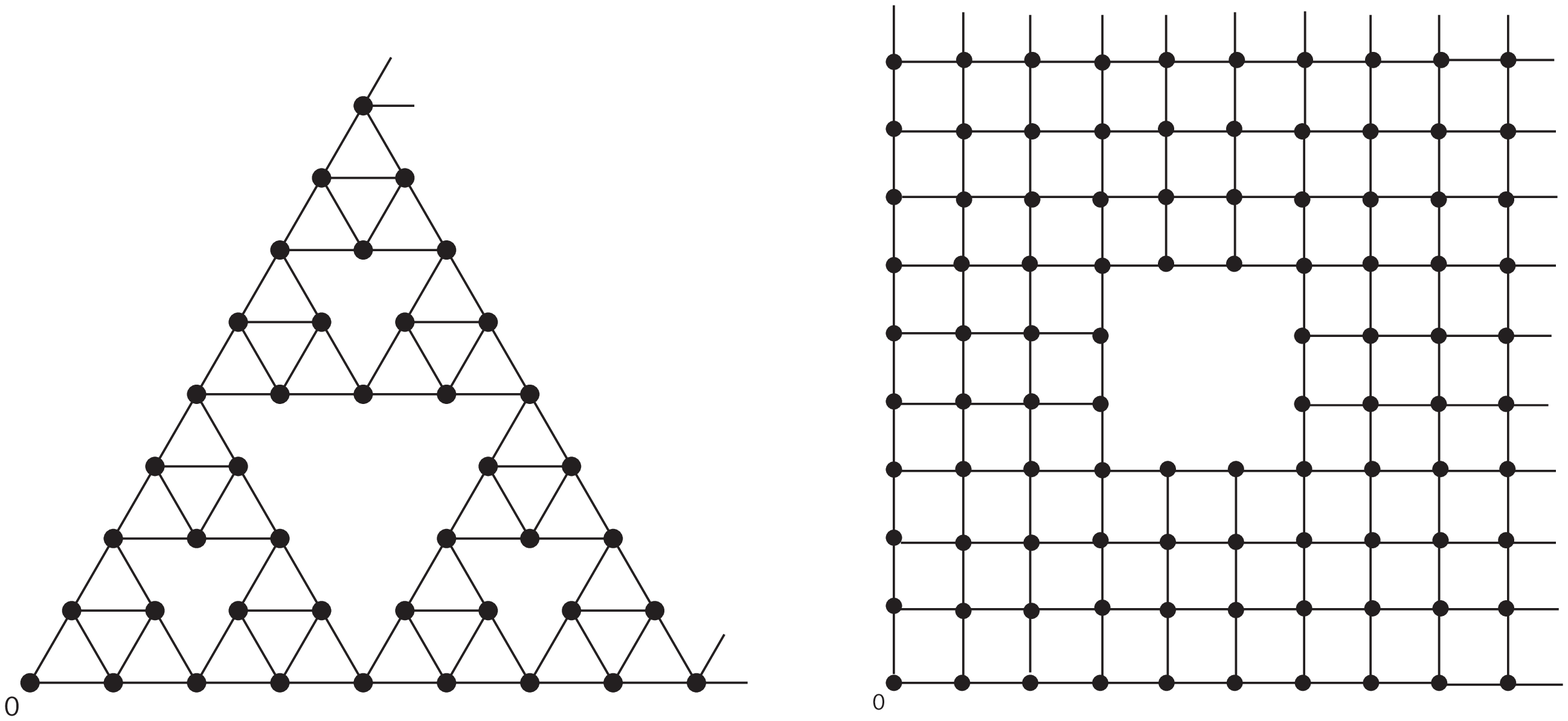, height=2in}}
\vspace{-.3cm}
\begin{center} 
Fig 1: ~2-dimensional Sierpinski gasket graph and carpet graph 
\end{center}
\subsubsection{Nested fractal graphs} 
Let $\Xi$ be the set of fixed points of $\{\psi_i\}_{i=1}^K$, and define 
\[V_0:=\left\{x\in \Xi:\mbox{ $\exists i,j\in\{1,\dots,K\}$, $i\neq j$ and $y\in \Xi$ such that $\psi_i(x)=\psi_j(y)$}\right\}.\]
Assume that $\# V_0\geq 2$ and set $\psi_{i_1\dots i_n}:=\psi_{i_1}\circ\dots\circ\psi_{i_n}$. $F$ is then called a nested fractal if the following holds. 
\begin{itemize}
  \item (Nesting) If $i_1\dots i_n$ and $j_1\dots j_n$ are distinct sequences in $\{1,\dots,K\}$, then
  \[\psi_{i_1\dots i_n}(F)\cap\psi_{j_1\dots j_n}(F)=\psi_{i_1\dots i_n}(V_0)\cap\psi_{j_1\dots j_n}(V_0).\]
  \item (Symmetry) If $x,y\in V_0$, then the reflection in the hyperplane $H_{xy}:=\{z\in \mathbb{R}^d:\:|z-x|=|z-y|\}$ maps 
  $\bigcup_{i_1,\dots,i_n=1}^K \psi_{i_1\dots i_n}(V_0)$ to itself.
\end{itemize}
We assume without loss of generality that $\psi_1(x)=L^{-1}x$ and that
the origin belongs to $V_0$. 
Let 
\bee\label{eq:vgnlwd}
V(G^N):=\bigcup_{i_1,\dots ,i_N=1}^{K}L^N\psi_{i_1\dots i_N}(V_0),~~~G:=\bigcup_{N=1}^{\infty}V(G^N).\eee
Next, define $B_0:=\{\{x,y\}:x\ne y\in
V_0\}$.  Then inside each $L^{N}\psi_{i_{1}\dots i_{N}}(V_0)$,  
$N\ge 0, 1\le i_1,\cdots,i_N\le K$, we place a copy of $B_0$ and denote by 
$B$ the set of all the edges determined in this way. 
Next, we assign $\mu_{xy}
=\mu_{yx}>0$ for each $\{x,y\}\in B$ in such a way that there exist $c_1,c_2>0$ such that
\[c_1\le \mu_{xy}=\mu_{yx}\le c_2,\qquad \forall \{x,y\}\in B.\]
We call the graph $(G,\mu)$ a nested fractal graph. 
A typical example is the 2-dimensional Sierpinski gasket graph 
in Fig 1 (where $L=2$). 
Let $d(\cdot,\cdot)$ be the graph distance on $G$,
$\{w_k\}_k$ the Markov chain for $(X,\mu)$, and
define the heat kernel as $p_k(x,y)=P^x(w_k=y)/\mu_y$.   
(Note that we consider the discrete time Markov chain here
in order to apply the results in \cite{BCK} to derive the resistance estimates \eqref{eq:res1}. Indeed,  
\eqref{eq:res1} can be obtained through both discrete and continuous time Markov chains.) 
It is known (see \cite{HK04} (also \cite{KumaHK} for the continuous setting)) that 
there exist constants $c_{3},\dots,c_{6}$ such that for
all $x,y\in G, k>0$
\bee\label{eq:HKup1} p_k(x,y) \leq c_{3} k^{-d_f/d_w} \exp\left(-c_{4}
\left(\frac{d(x,y)^{d_w}}{k}\right)^{1/(d_w-1)}\right), \eee
and for $k>d(x,y)$,
\bee\label{eq:HKlow1}
p_k(x,y)+p_{k+1}(x,y) \geq c_{5} k^{-d_f/d_w} \exp\left(-c_{6}
\left(\frac{d(x,y)^{d_w}}{k}\right)^{1/(d_w-1)}\right), \eee
where $d_w=\log (\rho K)/\log(L\eta)$, $d_f=\log K/\log(L\eta)$ with some constants $\rho>1$, $\eta\ge 1$. 
$d_f$ is called the \textit{Hausdorff dimension}
and $d_w$ is called the \textit{ walk dimension}.   
For the 2-dimensional Sierpinski gasket graph, $L=2, \eta=1$, $K=3$ and $\rho=5/3$.  
Noting that $d_w>d_f$ and that 
\bee\label{eq:vol1}
c_7R^{d_f}\le \mu (B(x,R))\le c_8 R^{d_f}, ~~~\forall x\in G, R\ge 1,
\eee
\eqref{eq:HKup1}, \eqref{eq:HKlow1} implies (see \cite[Theorem 1.3, Lemma 2.4]{BCK}) 
\bee\label{eq:res1}
R(x,y)\le c_9d(x,y)^{d_w-d_f},~~~
R(x,B^c(x,R))\ge c_{10}R^{d_w-d_f},~~\forall x,y\in G, ~\forall R\ge 1.
\eee
We now define a sequence of graphs $\{G^N\}_{N\geq 0}$ by setting $V(G^N)$ as above 
and $E(G^N):=\{\{x,y\}\in B: x,y\in V(G^N)\}$. Let $d_{G^N}(\cdot,\cdot)$ be the graph distance on $G^N$; 
one can easily see that
$d(x,y)\le d_{G^N}(x,y)$ for $x,y\in G^N$. (Note that $|x-y|\asymp d_{G^N}(x,y)^{\log L/\log (L\eta)}$ 
for $x,y\in G^N$ (cf. \cite[Section 3]{KumaHK}) and $\log L/\log (L\eta)$ is called the chemical-distance exponent.)

Let $B^N:=L^NV_0$. Clearly $R_{B^N}(x,y)\le R(x,y)$ for $x,y\in G^N$ and $d_{max}^N\asymp d_{G^N}(0, B^N)\asymp (L\eta)^N$. 
So \eqref{eq:res1} implies Assumption \ref{thm:assump1} (i),(ii) with $h(s)=s^{d_w-d_f}$, and 
\eqref{eq:vol1} with the self-similarity of the graph imply Assumption \ref{thm:assump1} (iii) with $\alpha=d_f$. 
We note that we can actually take $B^N$ arbitrary as long as $d_{max}^N\asymp (L\eta)^N$. \\

\subsubsection{Strongly recurrent Sierpinski carpet graphs}
Let $H_0=[0,1]^d$, and let $L \in \mathbb N$, $L\geq 2$ be fixed. 
Set ${\mathcal Q}=\{\Pi_{i=1}^d [(k_i-1)/L,k_i/L]: 1\le k_i\le L~
(1\le i\le d)\}$, let 
$L\le K \le L^d$ and let $\{\psi_i\}_{i=1}^K$
be a family of $L$-similitudes of
$H_0$ onto some element of ${\mathcal Q}$. We assume that the sets $\psi_i(H_0)$ are distinct, and as before
assume $\psi_1(x)=L^{-1}x$.  
Set $H_1=\cup_{i=1}^K\psi_i (H_0)$. Then, there exists a unique
non-void compact set $F \subset H_0$ such that $F=\cup_{i=1}^K\psi_i(F)$.
We assume $F$ is connected. $F$ is called a 
(generalized) Sierpinski carpet if the following hold (cf. \cite{BBconcar}): \\
(SC1) (Symmetry) $H_1$ is preserved by 
all the isometries of the unit cube $H_0$.\\
(SC2) (Non-diagonality) Let $B$ be a cube in $H_0$ which is the union of
$2^d$ distinct elements of ${\mathcal Q}$. (So $B$ has side length $2L^{-1}$.) Then
if $\mbox{Int} (H_1 \cap B)$ is non-empty, it is connected.\\
(SC3) (Borders included) $H_1$ contains the line segment $\{x: 0\leq x_1
\leq 1, x_2=\cdots =x_d=0 \}$.

The main difference from nested fractals 
is that Sierpinski carpets are infinitely ramified, i.e. 
$F$ cannot be disconnected by removing a finite number of points.

Let $V_0$ be a set of vertices in $H_0$ and define $V(G^N)$ and $G$ as in
(\ref{eq:vgnlwd}). Set $B_0:=\{\{x,y\}:x\ne y\in V_0, |x-y|=1\}$, and define 
$B$ and $\mu_{xy}$ as in the case of nested fractal graphs. 
We call the graph $(G,\mu)$ a Sierpinski carpet graph. A typical example is the 2-dimensional Sierpinski carpet graph in Fig 1.

It is known, 
see \cite{BB99} and also \cite{BBconcar} for the continuous setting, that 
\eqref{eq:HKup1}, \eqref{eq:HKlow1} hold, 
where $d_w=\log (\rho K)/\log L$, $d_f=\log K/\log L$ with some constant $\rho>0$. 
For the 2-dimensional Sierpinski gasket graph, $L=3$, $K=8$ and $\rho>1$. 
Let us restrict ourselves to the case $\rho>1$, namely $d_w>d_f$. 
In this case, since \eqref{eq:vol1} holds, we can show that \eqref{eq:HKup1} and \eqref{eq:HKlow1} imply \eqref{eq:res1} as before.
Arguing further as before, we have Assumption \ref{thm:assump1} (i)--(iii) with $h(s)=s^{d_w-d_f}$ and 
 $\alpha=d_f$.
         \subsection{Homogeneous random Sierpinski carpet graphs}
Let $\ell\ge 2$ and $I:=\{1,\cdots,\ell\}$. For each $k\in I$, let $\{\psi^k_i\}_{i=1}^{K_k}$ be a family of $L_k$-similitudes 
as in the definition of the Sierpinski carpet graphs. 
As before, we assume $\psi_1^k(x)=L_k^{-1}x$. For $\xi=(k_1,\cdots,k_n,\cdots)\in I^\infty$ and $n\in \mathbb N$, write $\xi|_N=(k_1,\cdots,k_N)\in I^N$, and let 
\bee\label{eq:vgnlwd22}
V(G_{\xi|_N}^N):=\bigcup_{{i_j\in \{1,\cdots,K_{k_j}\},}\atop{1\le j\le N}}L_{k_1}\cdots L_{k_N}
\psi^{k_N}_{i_N}\circ\dots\circ\psi^{k_1}_{i_1}(V_0),
~~~G_\xi:=\bigcup_{N=1}^{\infty}V(G_{\xi|_N}^N).\eee
Let $B_0:=\{\{x,y\}:x\ne y\in V_0, |x-y|=1\}$, and define 
$B=B_{\xi}$ as in the cases of nested fractal graphs and carpet graphs. For simplicity,
put weight $\mu_{xy}\equiv 1$ for each $\{x,y\}\in B$. 
We call the graph $(G_\xi,\mu_\xi)$ a homogeneous (random) Sierpinski carpet graph.

Fix $n\in \mathbb N$, $\xi|_n=(k_1,\cdots,k_n)\in I^n$, and let $B_n=L_{k_1}\cdots L_{k_n}$, 
$M_n=K_{k_1}\cdots K_{k_n}$.  
We write $R_n$ for the effective resistance between $\{0\}\times [0,B_n]^{d-1}\cap G_{\xi|_n}^n$ and 
$\{B_n\}\times [0,B_n]^{d-1}\cap G_{\xi|_n}^n$ in $G_{\xi|_n}^n$, and define $T_n=R_nM_n$. 
 Now set
\[
 d_f(n)= \frac{\log M_n}{\log B_n}, \qquad
   d_w(n) = \frac{\log{T_n}}{\log{B_n}}. 
\]
For $x\in G_\xi$
and $r\ge1$, let $V_d(x,r)$ be the number of vertices in the ball of radius $r$ centered at $x$ w.r.t. the graph
distance. It can be easily seen that 
\bee
c_1r^{d_f(n)}\le V_d(x,r)\le c_2 r^{d_f(n)}~~~~\mbox{if }~~B_n\le r< B_{n+1}, ~x\in G_{\xi}.\label{eq:volrand}\eee
         Define a time scale function $\tau: [1,\infty)\to [1,\infty)$ and resistance scale factor $h: [1,\infty)\to [1,\infty)$
as
\[\tau(s)=s^{d_w(n)},~ h(s)=s^{d_w(n)-d_f(n)}~~~~\mbox{if }~~T_n\le s< T_{n+1}.\]
We set $\tau(0)=h(0)=0$. Note that $\tau$ and $h$ satisfy the property in \eqref{eq:upvol} since $\ell<\infty$.  

Given these, it is possible to obtain heat kernel estimates similar to those in Theorem 6.3 and 
Lemma 6.7 of \cite{HKKZ} by tracking the proof in \cite{HKKZ} faithfully
(see the Appendix for a sketch).  
By making additional computations (similar to those in \cite[Lemma 3.19]{GT3}) in the proof of \cite[Lemma 3.10]{HKKZ}, 
we can obtain the following heat kernel estimates (cf. Remark after Theorem 24.6 in \cite{Kig}):  
There exist $c_{3},\cdots, c_{6}>0$ such that if $k\in \mathbb N$, $x,y \in G_\xi$, then   
\begin{eqnarray}
 p_k(x,y) &\leq &\frac{c_3}{V_d(x,\tau^{-1}(k))}\exp\Bigl(-c_4\big(\frac{\tau(d(x,y))}{k}\big)^{1/(\beta_1-1)}
\Bigr),\label{eq:ubt2}\\
p_k(x,y)+p_{k+1}(x,y) &\geq &\frac{c_5}{V_d(x,\tau^{-1}(k))}~~~\mbox{for }~k\ge c_6\tau(d(x,y)).\label{eq:lbt2}
\end{eqnarray}
Now assume the following limits exist and the inequality holds.
\bee d_f:=\lim_{n\to\infty}d_f(n),~~d_w:=\lim_{n\to\infty}d_w(n),~~d_w>d_f. \label{eq:strscho} 
\eee
Under this assumption, we have
\bee
c_7\frac{\tau(d(x,y))}{V_d(x,d(x,y))}\le R(x,y)\le c_8\frac{\tau(d(x,y))}{V_d(x,d(x,y))},~~~\forall x,y\in G_\xi.\label{eq:resfbue} 
\eee
The equivalence of \eqref{eq:ubt2}+\eqref{eq:lbt2} and \eqref{eq:resfbue} is proved in \cite{BCK} when $\tau(s)=s^\beta$ for
some $\beta\ge 2$ under some volume growth condition referred as $(VG(\beta_-))$. Here we need a generalized version of this under 
the doubling property of $\tau$. In fact, we only need \eqref{eq:ubt2}+\eqref{eq:lbt2} $\Rightarrow$ \eqref{eq:resfbue}, 
and the generalization of this direction is easy. 
Indeed, using \eqref{eq:ubt2} and \eqref{eq:lbt2}, we can obtain the scaled Poincar\'e inequality  and the lower bound of 
\eqref{eq:resfbue} similarly to the proof of \cite[Proposition 4.2]{BCK}
(with $\tau(s)$ replacing $s^\beta$ there). 
Under \eqref{eq:strscho}, a condition corresponding to $(VG((d_w)_-))$ in \cite{BCK} holds, so together with 
the scaled Poincar\'e inequality, we can obtain the upper bound of \eqref{eq:resfbue} similarly to the proof of 
\cite[Lemma 2.3 (b)]{BCK}. 

Now let $B^{\xi|_N}:=B_NV_0$. Clearly $R_{B^{\xi|_N}}(x,y)\le R(x,y)$ for $x,y\in G_{\xi|_N}^N$ and $d_{max}^N\asymp B_N$. 
So \eqref{eq:resfbue} implies Assumption \ref{thm:assump1} (i),(ii), and 
\eqref{eq:volrand}, \eqref{eq:strscho} with the homogeneity of the graph imply Assumption \ref{thm:assump1} (iii) with 
$\alpha=\max_nd_f(n)$. 
As before we can take $B^{\xi|_N}$ arbitrary as long as $d_{max}^N\asymp B_N$. 

Finally we will introduce randomness on this graph.
Let $(I^{\mathbb N}, {\cal F}, \mathbb P)$ be a Borel probability space where the measure $\mathbb P$ 
is stationary and ergodic for the shift operator 
$\theta: I^{\mathbb N}\to I^{\mathbb N}$ defined by $\theta((k_1,\cdots, k_n,\cdots))=(k_2,\cdots, k_n,\cdots)$. 
Then, by \cite[Proposition 7.1]{HKKZ} and the sub-additive ergodic theorem, one can prove the existence of the
first two limits in \eqref{eq:strscho}. 
Let $d_f^i, d_w^i$ be the Hausdorff dimension and the walk 
dimension for $G_{{\bf i}}$ where ${\bf i}=(i,i,i,\cdots)$ for $i\in I$. 
Let us consider a special case when $d=3, \ell=2$, and $\mathbb P$ is the Bernoulli probability measure
with $\mathbb P (\xi_1=1)=p$, $\mathbb P (\xi_1=2)=1-p$ for some $p\in [0,1]$.  
One can see that $d_f/d_w$ is a continuous function of $p$. 
Indeed, it can be easily seen that it is enough to prove $\lim_{n\to\infty}R_n/n$ is continuous 
for $p$. By the proof of \cite[Proposition 7.1]{HKKZ}, there exist $c_1,c_2>0$ such that we have 
\[\frac 1k \mathbb E\log (c_1R_k)\le \lim_{n\to\infty}\frac 1n R_n
\le \frac 1k \mathbb E\log (c_2R_k),~~~\mathbb P-\mbox{a.s.},
\]
for any $k\ge 1$ where $\mathbb E$ is the average over $\mathbb P$. 
Since $\mathbb E\log (c_iR_k)$, $i=1,2$ are continuous for $p$ (because the graph is finite), we 
obtain the desired continuity of $\lim_{n\to\infty}R_n/n$. 
So, when we choose the two carpets in such a way 
that $d_w^1>d_f^1$ and $d_w^2<d_f^2$ (which is possible, see \cite[Section 9]{BBconcar}), we are able to construct a one
parameter family of homogeneous random Sierpinski carpet graphs where $d_f/d_w$ is $\mathbb P$-a.e. an arbitrary fixed 
number between $d^1_f/d^1_w$ and $d^2_f/d^2_w$. In particular, there exists $p_*\in (0,1)$ such that \eqref{eq:strscho} holds
$\mathbb P$-a.e. for all $p<p_*$.

\appendix
\section{Appendix: Heat kernel estimates for Markov chains on homogeneous random Sierpinski carpet graphs}\label{appen}
In this appendix, we will briefly sketch the proof of \eqref{eq:ubt2} and \eqref{eq:lbt2}.
The Markov chain we consider here is the discrete time Markov chain.  

Set $V_n:=V(G_{\xi|_n}^n)$. 
We first define the Dirichlet form as follows.  
\[
{\cal E}_n(f,g):=\sum_{{x,y\in V_n}\atop{\{x,y\}\in B}}(f(x)-f(y))(g(x)-g(y)),\qquad\forall f,g: V_n\to {\mathbb R}.
\]

Given two processes $Y^1,Y^2$, defined on the same state space, we define a coupling time 
of $Y^1$ and $Y^2$ as 
\[T_C(Y^1,Y^2)=\inf \{t\ge 0: Y^1_t=Y^2_t\}.\]
Let $m\le n$. We call sets of the form 
$L_{k_1}\cdots L_{k_n}
\psi^{k_{n-m}}_{i_{n-m}}\circ\dots\circ\psi^{k_1}_{i_1}([0,1]^d)\cap V_n$ $m$-{\sl complexes}.
For $A\subset G_\xi$, define 
\begin{eqnarray*}
D^0_m(A)&=&\{\mbox{$m$-complex which contains } A\},\\
D^1_m(A)&=&D^0_m(A)\cup \{B:B \mbox{ is a $m$-complex}, D_m^0(A)\cap 
B \neq \emptyset\}. 
\end{eqnarray*}
Let $S_{B}^z$ denote the exit time from the set $B$, when the
process is started from the point $z$. 
\begin{thm}(Coupling)\label{thm:couple}
There exist $0<p_0<1$ and $K_0\in {\mathbb N}$ such that 
for each $x,y\in G_\xi$, there exist Markov chains $\bar w_t^x, \bar w_t^y$ with   
$\bar w_0^x=x$, $\bar w_0^y=y$ on $G_\xi$ whose laws are equal to the simple random walk that
satisfy the following: For $n> K_0$ and $y \in D^0_{n-K_0}(x)$,
\[ P(T_C(\bar w_t^x,\bar w_t^y)<\min\{S^x_{D^1_n(x)},S^y_{D^1_n(x)}\}) > p_0. \]
\end{thm}
The proof of the theorem follows in the same way as \cite[Section 3]{BBconcar},
as $G_\xi$ and $\bar w_t^x$ have enough symmetries for the argument there to work.

Once we have the coupling estimate, we can deduce the uniform (elliptic) Harnack inequality
as in \cite[Section 4]{BBconcar}. Let ${\cal L}$ be the infinitesimal 
generator associated with the simple random walk. 
\begin{thm}\label{thm:ehiuni}
There exists $c_1>0$ such that for each $x_0\in G_\xi$, and each $f: B(x_0,2R)\to [0,\infty)$
with ${\cal L} f(x)=0$ for all $x\in B(x_0,2R)= 0$, $R\ge 1$, it holds that 
\begin{equation}
\max_{x\in B(x_0,R)}f(x)\le c_1\min_{x\in B(x_0,R)}f(x).\label{eq:lhar}
\end{equation}
\end{thm}
\noindent
We next introduce the following Poincar\'e constant:
\begin{eqnarray*}
\lambda_n&=& \sup \{\sum_{x\in V_n}(u(x)-\langle 
u\rangle_{V_n})^2\;|\;  u: V_n\to \bR,\,
{\cal E}_n(u,u)=1\},
\end{eqnarray*}
where $\langle u\rangle_A= (\sharp A)^{-1} \sum_{x\in A}u(x)$
for any finite set $A$ and $u: V_n\to \bR$. 

The following proposition can be proved similarly to Proposition 3.1, Corollary 3.7 
of \cite{HKKZ} and (2.3), (4.4) of \cite{kz}. (Note that Theorem \ref{thm:ehiuni} is needed in the proof of \eqref{eq:B_re22}.)  
\begin{propn}
There exist constants $c_{1},\cdots, c_{4}>0$ such that for each $n,m\in {\bf N}$,
\begin{eqnarray}
c_{1} R_nR_{\theta^n\xi\vert m}&\le &R_{n+m}\le
c_{2}R_nR_{\theta^n\xi\vert m}, \label{eq:B_re}\\
c_3\lam_n&\le & T_n\le c_4\lam_n.\label{eq:B_re22}
\end{eqnarray}
\label{thm:limre}
\end{propn}

\begin{lem} \label{lem:updiag}
There is a constant $c_{}$ such that if 
$T_{n-1} \leq t \leq T_{n}$, then
\begin{equation}
  p_t(x,y)\leq c M_n^{-1}.
\end{equation}
\end{lem}
\noindent {\it Proof.}  From the definition of the Dirichlet form and the Poincar\'e constant,
the proof is similar to \cite[Theorem 3.3]{kz} by using Proposition \ref{thm:limre}. \qed 

The next lemma can be proved similarly to \cite[Lemma 3.8]{HKKZ}. 
\begin{lem}\label{thm:averes}
There exist $c_1,c_2>0$ such that 
\[
c_{1} T_r \leq ES^z_{D^1_r(x)}~\mbox{for all } z\in D^0_{r}(x),~~~~
ES^z_{D^1_r(x)} \leq c_{2} T_r~\mbox{for all } z\in D^1_{r}(x).
\]
\end{lem}
\noindent Since $S_{D^1_l(x)} \leq t + 1_{(S_{D^1_l} >t)}(S_{D^1_l}-t)$ 
we have, from Lemma \ref{thm:averes},
\[
c_1T_l\le ES^z_{D^1_l} \leq t + E \bigl[ 1_{(S^z_{D^1_l}>t)} E[S^{X_t}_{D^1_l}] 
 \bigr]   \leq t + P(S^z_{D^1_l}>t) c_2 T_l~~~\mbox{for } t \geq 0, z\in D_l^0(x). 
 \]
Thus, we deduce the following: there exist $c_3>0$, 
$c_4\in (0,1)$ such that
\begin{equation}
   P(S^z_{D^1_l(x)} \leq t) \leq c_3 T_l^{-1} t  + c_4~~\mbox{for } t \geq 0,
   z\in D_l^0(x).
\label{eq:ubwn}
\end{equation}
We can improve this to an exponential estimate on 
$P(S^z_{D^1_l(x)} \leq t)$. In order to do this we define 
the following function of time and space,
\begin{equation}\label{eq-A6}
k=k(n,l) = \inf\bigl\{l'\le n :\frac{T_{l'}}{B_{l'}} 
   \geq \frac{T_l}{B_n} \bigr\}. 
\end{equation}

The next lemma corresponds to \cite[Lemma 3.10]{HKKZ}. 
Since the labeling here differs from that in \cite{HKKZ}, we give the proof. 
\begin{lem}  \label{lem:chainhit}
There exist constants $c_{1}, c_{2}$ such that if $k=k(n,l)$ as in \eqref{eq-A6}
then for
all $x\in E$, and $n,l \ge 0$,
\begin{equation}
 P(S^x_{D^1_n(x)} \leq T_l) \leq c_{1} \exp{(-c_{2}B_{n}/B_k)}.
  \label{eq:expest1}
\end{equation}
\end{lem}
\noindent {\it Proof.}
If $l'\le n$, then for the simple random walk to cross one $n$-complex it must
cross at least $N=B_{n}/B_{l'}$, $l'$-complexes. So, there exists $0<c<1$ such 
that 
\[ S^x_{D^1_n(x)} \geq  \sum_{i=1}^{cB_{n}/B_{l'}} V_i^{x_i}, \]
where $x_i$ depend only on $V_1^{x_1},\ldots,V_{i-1}^{x_{i-1}}$, and $V_i^y$ have the same distribution as
$S^y_{D^1_{l'}(y)}$. The deviation estimate \cite[Lemma 1.1]{baba} states that if $P(V_i^y< s) \leq p_0 + \alpha s$, where 
$p_0\in (0,1)$ and $\alpha>0$, then
\begin{equation}
  \log P \bigl( \sum_1^{cN} V_i^{x_i} \leq t) 
   \leq 2(\alpha c_1Nt/p_0)^{1/2} - c_2N \log(1/p_0).
\label{eq:l11}
\end{equation} 
Thus, using (\ref{eq:ubwn}) and (\ref{eq:l11}), we have 
\begin{equation}
\label{eq:ubwn2}
 \log P(S^x_{D^1_n(x)} \leq T_l) \leq
c_{3} (B_{n}/B_{l'})^{1/2}[ (T_{l}/T_{l'})^{1/2} -c_{4} (B_{n}/B_{l'})^{1/2} ].
\end{equation}
Given $k=k(n,l)$ as above, there exists $c_{5}$ and $k_0$ such that 
$k \leq k_0 \leq k+c_5$, and
\[ (T_{l}/T_{k_0})^{1/2} < \frac 12 c_{4} (B_{n}/B_{k_0})^{1/2} . \]
Provided $k_0 \le n$ we deduce
\[ \log P(S^x_{D^1_n(x)}\leq T_l) \leq -\frac 12 c_{3} c_{4} 
B_{n}/B_{k_0}. \]
Choosing $c_6$ large enough we have 
 $1 < c_{6} \exp(-c_{2} B_{n}/B_k)$ whenever 
$k > n-c_5$, so that (\ref{eq:expest1}) holds in all cases. \qed

\begin{thm} \label{thm:mainub}
There exist constants $c_{1},c_{2}$ such that if $k\in {\mathbb N}$, $x,y \in G_\xi$, 
and $n,m$ satisfy
\begin{equation}
 T_{n-1} \leq t < T_{n}, \quad B_{m-1} \leq d(x,y) < B_{m},
\label{eq:conxyt}
\end{equation}
and $k=k(m,n)$, then
\begin{equation}
 p_t(x,y) \leq c_{1} t^{-d_f(n)/d_w(n)} \exp \Bigl( 
    -c_{2} \bigl(\frac{d(x,y)^{d_w(k)}}{t} \bigr)^{1/(d_w(k)-1)}
\Bigr).
\label{eq:ubt}\end{equation}
\end{thm}

\noindent {\it Proof.} 
Noting that $M_n^{-1} \leq c t^{-d_f(n)/d_w(n)}$, this 
is proved from Lemma \ref{lem:updiag} and Lemma \ref{lem:chainhit} 
by the same argument as in Theorem $6.9$ of \cite{BBconcar}. \qed

Note that the bound (\ref{eq:ubt}) may also be written in the form 
\[
 p_t(x,y) \leq c M_n^{-1} \exp(-c^\prime B_{n} /B_k),
\]
 where $m,n$ satisfy (\ref{eq:conxyt}), and $k=k(m,n)$ as in \eqref{eq-A6}.  The upper bound 
\eqref{eq:ubt2} can be obtained from this using \eqref{eq:volrand}. \\

The lower bound is obtained in the following procedure. 
\begin{lem} \label{lem:diag}
There exists a constant $c_{1}$ such that if $T_n\leq t $
then 
 \begin{equation}
p_t(x,x) \geq c_{1} M_n^{-1} \quad \textrm{ for all } x\in G_\xi. 
\label{eq:onlub1}
\end{equation}
\end{lem}
\noindent {\it Proof.} Using Lemma \ref{lem:updiag} and \eqref{eq:volrand},
a standard argument gives the desired estimate. See for instance 
 \cite[Lemma 5.1]{bah}. \qed 
 \begin{lem} \label{lem:ndiag}
There exist $c_{1}, c_{2}$ such that if $T_{n-1} <t\le T_{n}$, then
\[
   p_t(x,y) \geq c_{1} M_n^{-1} \quad 
   \textrm{ whenever } d(x,y) \leq c_{2} B_n.
\]
\end{lem}
\noindent {\it Proof.} Using Theorem \ref{thm:ehiuni} and Lemma \ref{lem:diag},
this can be proved similarly to the proof of \cite[Proposition 6.4]{BB99}. \qed 

\noindent We can deduce \eqref{eq:lbt2} from this and  \eqref{eq:volrand}.
\qed

\bigskip
\noindent {\bf Acknowledgment. } The authors thank 
D. Croydon for valuable comments.

\vskip 0.3truein

\noindent {\bf Takashi Kumagai}

\noindent Research Institute for Mathematical Sciences,
Kyoto University, Kyoto 606-8502, Japan.

\noindent E-mail: {\tt kumagai@kurims.kyoto-u.ac.jp}

\bigskip
\noindent {\bf Ofer Zeitouni}
\noindent University of Minnesota, School of Mathematics, 
206 Church Street SE, Minneapolis, MN 55455, USA.

and

\noindent Weizmann Institute of Science, Faculty of Mathematics,
POB 26, Rehovot 76100, Israel.

\noindent E-mail: {\tt zeitouni@math.umn.edu}

\begin{thebibliography}{99}
             \bibitem{Adl}R.J. Adler.\newblock {\em An introduction to continuity, extrema, and related topics for general Gaussian processes.}
               \newblock IMS Lecture Notes-Monograph Series, 1990.
               \bibitem{baba}
M.T.~Barlow and R.F.~Bass.
\newblock Construction of Brownian motion on the Sierpinski carpet.
\newblock {\em Ann. Inst. H. Poincar\'e}, {\bf 25}, 225--257, 1989.
             \bibitem{BB99} M.T. Barlow and R.F. Bass. 
               \newblock {\em Random walks on graphical Sierpinski carpets.} 
               \newblock In: Random walks and discrete potential theory, pp. 26-55, Cambridge, Cambridge Univ. Press, 1999.
             \bibitem{BBconcar} M.T. Barlow and R.F. Bass.  Brownian motion 
               and harmonic analysis on Sierpinski carpets. 
               {\em Canadian J. of Math.} {\bf 51} (1999), 673--744.
             \bibitem{BCK} M.T. Barlow, T. Coulhon and T. Kumagai.
               \newblock Characterization of sub-Gaussian heat kernel estimates on strongly recurrent graphs.
               \newblock {\em Comm. Pure Appl. Math., \bf 58} (2005), 1642--1677.              
\bibitem{bah}M.T. Barlow and B.M. Hambly.
\newblock Transition density estimates for Brownian motion on
scale irregular Sierpinski gasket
\newblock {\em Ann. Inst. H. Poincar\'e}, {\bf 33}, 531--557, 1997.
               \bibitem{BZrec}M. Bramson and O. Zeitouni. 
               \newblock Tightness for a family of recursion equations. 
               \newblock {\em Ann. Probab. \bf 37} (2009), 615--653.
             \bibitem{BZ10}M. Bramson and O. Zeitouni. 
               \newblock Tightness of the recentered maximum of the two-dimensional
               discrete Gaussian Free Field. 
               \newblock {\em Comm. Pure Appl. Math. \bf 65} (2012), 1--20.
             \bibitem{DLP} J. Ding, J. R. Lee and Y. Peres,
               Cover times, blanket times and majorizing measures.
                \newblock {\em Annals Math. \bf 175} (2012), 1409--1471. 
             \bibitem{EKMRS} N. Eisenbaum, H. Kaspi, M. B. Marcus, J. Rosen and
               Z. Shi. \newblock
               A Ray-Knight theorem for symmetric Markov processes.
               \newblock {\em Annals. Probab. \bf 28} (2000), 1781--1796.
             \bibitem{GT3} A. Grigor'yan and A. Telcs. 
               \newblock Two-sided estimates of heat kernels on metric measure spaces. \newblock {\em Ann. Probab. \bf 40} (2012), 1212--1284.
             \bibitem{HK04}B.M. Hambly and T. Kumagai.
               \newblock {\em Heat kernel estimates for symmetric random walks on a class of fractal graphs and stability under rough isometries.} 
               \newblock In: Fractal geometry and applications, Proc. of 
               Symposia in Pure Math. {\bf 72}, Part 2, pp. 233--260, Amer. Math. Soc. 2004.
             \bibitem{HKKZ}B.M. Hambly, T. Kumagai, S. Kusuoka and X.Y. Zhou.
               \newblock Transition density estimates for diffusion processes on homogeneous random Sierpinski carpets.
               \newblock {\em J. Math. Soc. Japan \bf 52} (2000), no. 2, 373--408. 
             \bibitem{Kig}J. Kigami. 
\newblock Resistance forms, quasisymmetric maps and heat kernel estimates. 
               \newblock {\em Memoirs Amer. Math. Soc. \bf 216}, no. 1015 (2012). \bibitem{TKSF}T. Kumagai. \newblock Random walks on disordered media and their scaling limits.\newblock St. Flour Lecture Notes (2010). 
               {\tt http://www.kurims.kyoto-u.ac.jp/\string~kumagai/StFlour-TK.pdf}.
             \bibitem{KumaHK}
T. Kumagai.  
\newblock Estimates of transition densities for Brownian motion on nested fractals. 
\newblock {\em Probab. Theory Relat. Fields \bf 96} (1993), 205--224. 
\bibitem{kz}
S.~Kusuoka and X.Y.~Zhou.
\newblock Dirichlet forms on fractals: Poincar\'e constant and resistance.
\newblock {\it Probab. Theory Relat. Fields}, {\bf 93}, 169--196, 1992. 
             \bibitem{LT}
M. Ledoux and M. Talagrand.
\newblock {\em Probability in Banach spaces.
Isoperimetry and processes.} \newblock
  Springer-Verlag, Berlin,  (1991).
\bibitem{LP}
\newblock {Lyons, R. {\rm with} Peres, Y.} (2012).
\newblock {\em Probability on Trees and Networks}.
\newblock Cambridge University Press.
\newblock In preparation. Current
  version available at \hfil\break
  {\tt http://mypage.iu.edu/\string~rdlyons/}.
  
               \bibitem{She}S. Sheffield. 
               \newblock Gaussian free fields for mathematicians. 
               \newblock {\em Probab. Theory Rel. Fields. \bf 139} (2007),
               521--541.
               \end{thebibliography}
\end{document}